\documentclass[leqno, 11pt]{amsart}
\numberwithin{equation}{section}

\usepackage{amsmath,amsthm,amsfonts,amssymb}
\usepackage[backref=page]{hyperref}
\hypersetup{colorlinks=true,citecolor=blue,linkcolor=blue,urlcolor=blue,pdfstartview=FitH,
pdftitle=Sums of Hecke eigenvalues over  quadratic polynomials}

\theoremstyle{plain}
\newtheorem{satz}{Theorem}
\newtheorem{lemma}{Lemma}
\newtheorem{prop}[lemma]{Proposition}

\theoremstyle{definition}

\DeclareMathOperator{\sgn}{sgn}

\renewcommand{\geq}{\geqslant}
\renewcommand{\leq}{\leqslant}

\newcommand{\abs}[1]{\lvert#1\rvert}

\begin{document}

\title{Sums of Hecke eigenvalues over quadratic polynomials}

\author{Valentin Blomer}
\address{Department of Mathematics, University of Toronto, 40 St. George Street, Toronto, Ontario,
 Canada, M5S 2E4} \email{vblomer@math.toronto.edu}
\thanks{Author supported by NSERC grant 311664-05}


\begin{abstract} Let $f(z) = \sum_n a(n) n^{(k-1)/2} e(nz) \in S_k(N, \chi)$ be a cusp form for $\Gamma_0(N)$, weight $k \geq 4$ and character $\chi$. Let $q(x) = x^2+sx+t\in \Bbb{Z}[x]$ be a quadratic polynomial. It is shown that
\begin{displaymath}
  \sum_{n \leq X} a(q(n)) =  cX + O_{f, q, \varepsilon}(X^{6/7+\varepsilon})
\end{displaymath} 
for some constant $c = c(f, q)$. The constant vanishes in many (but not all) cases, for example if $k$ is even or if $\Delta = s^2-4t > 0$. On the way a Kuznetsov formula for half-integral weight and entries having different sign is derived.
 \end{abstract}


\maketitle


\section{Introduction}

While most classical arithmetic functions are reasonably well understood on average, the situation becomes much harder if one considers sums over sparse sequences, for example values of a polynomial $q \in \Bbb{Z}[x]$ with $\text{deg } q \geq 2$. One of the most challenging and famous problems in this direction is the asymptotic evaluation of $\sum_{n \leq X} \Lambda(n^2+1)$. For the divisor function, Hooley \cite{Ho} showed 
\begin{equation}
  \sum_{n \leq X} \tau(n^2+a) = c_1(a) X\log X + c_2(a)X + O(X^{8/9} (\log X)^3) 
\end{equation}
for any fixed $a \in \Bbb{Z}$ such that $-a$ is not a perfect square (with the convention $\tau(m) = 0$ if $m \leq 0$). He remarks that a refinement of the method can be applied to evaluating $\sum_{n \leq X} \tau(an^2 + bn + c)$ and also to $\sum_{n \leq X} r(an^2 + bn + c)$ where $r(m)$ is the number of representations of $m$ as a sum of two squares. The error term in (1.1) was improved by Bykovski\u i \cite{By} to $O(X^{2/3+\varepsilon})$. Nothing of type (1.1) is known for cubic or higher degree polynomials (cf. \cite{Ki}). If one replaces $\tau$ by $\tau_3$, Friedlander and Iwaniec (\cite{FI}, see also \cite{T}) have recently proved that
\begin{displaymath}
  \sum_{\substack{n^2 + m^6 \leq X\\ (n, m) = 1}} \tau_3(n^2 + m^6) = c X^{2/3} (\log X)^2 + O(X^{2/3}(\log X)^{7/4+\varepsilon}).
\end{displaymath}\\

Let $N \in \Bbb{N}$, $k \geq 4$ an integer, and $\chi$ a Dirichlet  character mod $N$ satisfying $\chi(-1) = (-1)^k$.  Let $f \in S_k(N, \chi)$ be a holomorphic cusp form of weight $k$ and character $\chi$ for $\Gamma_0(N)$   (possibly, but not necessarily, an eigenform of the Hecke operators) with Fourier expansion
\begin{equation}
  f(z) = \sum_{n=1}^{\infty} a(n) n^{(k-1)/2} e(nz), \quad \Im z > 0.
\end{equation}
Let $q(x) = x^2+sx+t \in \Bbb{Z}[x]$ be an integral monic quadratic polynomial with discriminant $\Delta = s^2-4t$. In this paper we want to prove the following analogue of (1.1).

\begin{satz} With $f \in S_k(N, \chi)$ as {\rm (1.2)} and $q \in \Bbb{Z}[x]$ as above, one has
\begin{displaymath}
  \sum_{n \leq X} a(q(n))  = cX +O_{f, q, \varepsilon}(X^{6/7+\varepsilon})
\end{displaymath}
for any $\varepsilon > 0$ and $X \geq 1$ and some constant $c = c(f, q)\in \Bbb{C}$. We have $c=0$ if $k$ is even or $\Delta > 0$ (or both). 
\end{satz}

Although not directly related to the sums occurring in Theorem 1, it may be noted that sums of the type $\sum_{n=1}^{\sqrt{|D|}}\lambda_f(n^2-D)$ have recently been investigated in N.\ Templier's thesis in connection with Heegner points on elliptic curves.  Theorem 1 is an immediate corollary of  the following smoothed version. Let $\theta < 1/2$ be any lower bound for the Selberg eigenvalue conjecture (Ramanujan-Petersson conjecture at the archimedean place), i.e.\  for an eigenvalue $\lambda = 1/4 + t^2$ of the Laplacian $y^2(\partial_x^2 + \partial_y^2)$ on a modular surface $\Gamma_0(N)\backslash \Bbb{H}$ one has $|\Im t| \leq \theta$. By \cite{KS} we know that $\theta = 7/64$ is an admissible  constant. 

\begin{satz} Let $1 \leq P \leq X^{1/2}$, and let $w$ be a smooth function compactly supported on $[X/2, 2X]$ satisfying $w^{(j)} \ll_j (P/X)^j$ for all $j \geq 0$. Then with $f \in S_k(N, \chi)$  as in {\rm (1.2)} and $q\in \Bbb{Z}[x]$ as above, one has
 \begin{equation}
  \sum_n a(q(n)) w(n)     = c \int_0^{\infty} w(y) dy +   O_{f, q, \varepsilon}\left(\left(X^{1/2+\theta} P^{3/2-\theta} + X^{1/2}P^{5/2}\right)X^{\varepsilon}\right). 
\end{equation}
for some constant $c \in \Bbb{C}$ for any $\varepsilon > 0$.  We have $c=0$ if $k$ is even or $\Delta > 0$ (or both). 
\end{satz}

\textbf{Remarks:} 1) If $\Delta > 0$ and $k$ is odd, the constant $c$ may or may not vanish and has a complicated structure. We shall give an explicit expression for $c$ if  $f$ is a Poincar\'e series $P_m$ (cf.\ (2.1) below) with $m$ not a square: Let $\mathcal{B} = \{f_j(z) = \sum_n a_j(n) e(nz)\}$ be an orthonormal basis of  the space $M_{1/2}(4[N, 4], \chi\chi_{-4})$ of half-integral weight modular forms of weight $4[N, 4]$ and character $\chi\chi_{-4}$. Let $\delta(h) = 2$ if $h \not=0$ and $\delta(0) = 1$, and let $T_{\nu}(x) = \cos(\nu \arcsin(x))$ be the $\nu$-th Chebychev polynomial. Then if $f = P_m$ for $m$ not a square, $k$ is odd, and $\Delta < 0$, the constant $c$ is given by
\begin{equation}
  c = \frac{8\sqrt{2} i}{\pi^{1/4}(k-1)} \sum_{\substack{h \equiv s \, (2)\\ 0 \leq h < 2 \sqrt{m}}} \delta(h) \sum_{f_j \in \mathcal{B}} a_j(|\Delta|)\bar{a}_j(4m-h^2) T_{k-1}\left(\frac{h}{2\sqrt{m}}\right).
\end{equation}
 There are cases with $c \not= 0$, cf.\ the example after the proof of Theorem 2, where we consider the situation $f = P_3 \in S_5(12, \chi_{-4})$ and $q(x) = x^2+x+1$. The structure of $c$ is strongly reminiscent of certain trace formulas for Hecke operators. Probably the constant $c$ has a more intrinsic definition; in particular, there should be an expression for any given $f$, not necessarily a Poincar\'e series, which clarifies its meaning. It is not hard to see (cf.\ section 3) that any $f \in S_k(N, \chi)$ can be written as a finite linear combination of Poincar\'e series $P_m$ with $m$ not a square (possibly with a remainder that contribues negligibly to (1.3)), so in principle we can compute $c$ in any given situation.  We leave it as an interesting open question to shed more light on its underlying structure. It is also interesting to note that the sums occurring in Theorems 1 and 2 are somehow able to detect the exceptional eigenvalue $\lambda = i/4$ of the half-integral weight Laplacian.\\

2) It should  be noted that there are no restrictions on the quadratic polynomial $q$, for example $q$ does not need to be irreducible over $\Bbb{Q}$. In particular, applying Theorem 2 to Hecke eigenforms with $q(x) = x(x+h)$, $h > 0$ fixed (so $\Delta = h^2 > 0$), and a weight function $w$ as above, we recover easily  best-possible individual bounds for smooth shifted convolution sums
\begin{displaymath}
  \sum_{n} a(n) a(n+h) w(n) \ll_{h, f, P, \varepsilon} X^{1/2 + \theta +\varepsilon}.
\end{displaymath}



3) It is possible to include the case $k=3$ with slightly more careful estimations. The interesting case $k=2$ would require much more effort, since  we would obtain several not convergent or not absolutely convergent series in the course of the proof. It is also possible to treat more general polynomials $q(x) = rx^2 + sx +t \in \Bbb{Z}[x]$ with $r > 0$ by the same method. With considerably more work it might be possible to improve the $P$-exponent of the second term on  the right-hand side of (1.3) somewhat; optimistically, one might be able to replace $P^{5/2}
$ with  $P^{2}$.\\ 

Let us briefly outline the method of the proof. In principle Theorem 1 should be of comparable difficulty to (1.1) in some respects, except, of course, that the Fourier coefficients $a(n)$ do not  allow a decomposition of the form $\tau= \textbf{1} \ast \textbf{1}$, which is the starting point of the argument in \cite{Ho} and \cite{By}. Thus Hooley's and Bykovski\u i's methods do not seem to be  applicable in our situation.  For holomorphic\footnote{Unfortunately, this trick does not work for Maa{\ss} forms.} cusp forms $f$, however, there is some kind of weak substitute of the decomposition $\tau =\textbf{1} \ast \textbf{1}$. We can write $f$ as a linear combination of Poincar\'e series and accordingly replace $a(q(n))$ with 
\begin{equation}
  \sum_{N \mid c} \frac{1}{c}S_{\chi}(q(n),m, c) J_{k-1}\left(\frac{\sqrt{q(n)m}}{c}\right)
\end{equation}
for certain numbers $m$. We can now evaluate the $n$-sum by Poisson summation. This is especially clean, if $c$ is a multiple of $4$,  which can always be arranged by embedding $S_k(N, \chi)$ into $S_k(4N, \chi)$, say. Summing (1.5) over $n$ with $w$ as in Theorem 2, we obtain something roughly of the form
\begin{equation}
  \sqrt{X} \sum_{4N \mid c} \frac{1}{c} K_{\chi}(m', -\Delta, c) \tilde{w}\left(\frac{c}{X}\right)
\end{equation}
for certain integers $m'$ where  $\tilde{w}$ has bounded support, and $K_{\chi}$ is a Kloosterman sum with theta multiplier. It is the key observation that Poisson summation together with quadratic reciprocity translates sums of Kloosterman sums over a quadratic polynomial into sums of Kloosterman sums with theta multiplier. We are now in a position to apply Kuznetsov's sum formula for half-integral weight to exploit cancellation in sums of type (1.6). Here two  problems can occur: one of the entries of the Kloosterman sum can vanish in which case the Kloosterman sum degenerates and Kuznetsov's formula is not available. These cases can either be excluded at the beginning of the argument or treated directly without appeal to Kuznetsov's formula. Secondly, the Laplacian of half-integral weight  has an exceptional eigenvalue $\lambda = 3/16$ coming from holomorphic modular forms of weight 1/2 and 3/2. Therefore the Kloosterman zeta-function $\sum_c K_{\chi}(m, n; c) c^{-2s}$ may have a pole at $s = 3/4$ that needs to be investigated carefully and may produce a main term. Using Shimura's correspondence and  the Kim-Sarnak bound \cite{KS}, the next eigenvalue is at least $\geq 1/4 - (\theta/2)^2$, so that the error term in (1.6) can be bounded by $X^{1/2 + \theta +\varepsilon}$. The $P$-dependence in (1.3) requires some careful analysis. \\

On the way,  we derive a Kuznetsov formula for half-integral weight where the two entries in the Kloosterman sum $K_{\chi}(a, -b; c)$ have opposite signs (Proposition 2).  Although in principle the method is well-known, this result does not seem to be in the literature. For example, with this formula at hand  one can easily extend a result of Bir\'o \cite[Theorem 2]{Bi}  and obtain an analogous formula for sums of a  Maa{\ss} form over generalized Heegner points in terms of Fourier coefficients of its Shimura lift. \\



\textbf{Acknowledgement.} I would like to thank Peter Sarnak and Kumar Murty for helpful discussions. 

\section{Notation and Lemmas}

\begin{lemma} [Poisson summation] Let $c \in \Bbb{N} $, $u \in \Bbb{Z}$, and $f$ be a Schwartz-class function. Then
\begin{equation}
  \sum_{\substack{n \in \Bbb{Z}\\ n \equiv u\, (c)}} f(n) = \frac{1}{c}\sum_{h \in \Bbb{Z}} \hat{f}\left(\frac{h}{c}\right) e\left(\frac{hu}{c}\right) 
  \end{equation}
where $\hat{f}(y) = \int_{\Bbb{R}} f(x) e(-xy)dx$ is the Fourier transform of $f$.
\end{lemma}

\textbf{Proof.} See e.g.\ \cite[(4.25)]{IK}.  \\

Let $N$ be a positive integer, $k \geq 4$ an integer, and let $\chi$ be a character mod $N$ satisfying $\chi(-1) = (-1)^k$.  
%
We define the Poincar\'e series
\begin{displaymath}
  P_m(z) = \sum_{\gamma =\left(\begin{smallmatrix}a &b\\ c &d\end{smallmatrix}\right) \in \Gamma_{\infty}\backslash \Gamma_0(N)} \bar{\chi}(\gamma) (cz+d)^{-k}e(m\gamma z), \quad m \in \Bbb{N}
\end{displaymath}
where $\Gamma_{\infty} = \left\{\left(\begin{smallmatrix} 1 & b \\ & 1\end{smallmatrix}\right) \mid b \in \Bbb{Z}\right\}$ is the stabilizer of the cusp $\infty$ and $\chi(\gamma) = \chi(d)$ if $\gamma =\left(\begin{smallmatrix}a &b\\ c &d\end{smallmatrix}\right)$.  Their Fourier expansion is given by (cf. \cite[Lemma 14.2]{IK})
\begin{displaymath}
 m^{(k-1)/2} P_m(z) = \sum_{n =1}^{\infty} \left(\delta_{nm}+2 \pi i^{-k}\sum_{N \mid c} \frac{1}{c} S_{\chi}(m, n; c) J_{k-1}\left(\frac{4\pi\sqrt{nm}}{c}\right)\right) n^{(k-1)/2} e(nz).
\end{displaymath}
Here
\begin{equation}
  J_{\nu}(x) = \sum_{j=0}^{\infty} \frac{(-1)^j}{j!\Gamma(j+1+\nu)} \left(\frac{x}{2}\right)^{2j+\nu}, \quad x \in \Bbb{C} \setminus (-\infty, 0]
\end{equation}
is the $J$-Bessel function and 
\begin{displaymath}
  S_{\chi}(m, n; c) := \left.\sum_{d \, (c)}\right.^{\ast} \chi(d) e\left(\frac{m d + n\bar{d}}{c}\right)
\end{displaymath}
is the twisted Kloosterman sum. As usual, the star indicates summation over residue classes coprime to the modulus. 
For $f(z) = \sum_n a(n) n^{(k-1)/2}e(nz) \in S_k(N, \chi)$ we have
\begin{equation}
  a(n) = \frac{(4\pi \sqrt{n})^{k-1}}{\Gamma(k-1)} \langle f, P_n \rangle.
\end{equation}
For $\nu \in \Bbb{Z}$, the $J_{\nu}$-Bessel function satisfies the bound
\begin{equation}
  J_{\nu}(x)  
    \ll_{\nu} \min(x^{|\nu|}, x^{-1/2}). 
\end{equation}
By convexity, this gives $J_{\nu}(x) \ll \min_{-1/2 \leq \ell \leq |\nu|} x^{\ell}$. We will often use this bound for various suitable (possibly non-integral) values of $\ell$.  For arbitrary $\nu \in \Bbb{C}$ one has
\begin{equation}
  J_{\nu}'(x) = \frac{1}{2}(J_{\nu-1}(x) + J_{\nu+1}(x)), 
\end{equation}
as well as the asymptotic expansion (\cite[8.451.1+7+8]{GR})
\begin{equation}
  J_{\nu}(x) = \frac{1}{\sqrt{x}} e\left(\frac{x}{2\pi}\right)C_{1}(x) +\frac{1}{\sqrt{x}} e\left(-\frac{x}{2\pi}\right)C_{2}(x) +O_{\nu,  A}(x^{-A}) 
\end{equation}
for any $A \geq 0$ and $\nu \in \Bbb{C}$, where $C_{1, 2}(x)$ are smooth functions satisfying $C_{1, 2}^{(j)}(x) \ll_{j,  \nu, A} x^{-j}$ for all $j \geq 0$. \\



An important ingredient for the results of this paper is Kuznetsov's sum formula. The exact statement needs some preparation. For the rest of the paper let  $\kappa \in  \{1/2, 3/2\}$.  Let $\Delta_{\kappa} :=  y^2(\partial_x^2 + \partial_y^2) -i\kappa y\partial_x$ be the Laplacian of weight $\kappa$. Let $\Lambda_{\kappa} := \kappa/2 + y(i\partial_x - \partial_y)$ be the Maa{\ss} lowering operator. Finally let $X$ be the reflection operator $(Xf)(z) = f(-\bar{z})$ for $f : \Bbb{H} \rightarrow \Bbb{C}$. Assume $4 \mid N$.  For odd $d$ let $\epsilon_d = 1$ if  $d \equiv 1$ (mod 4) and $\epsilon_d = i$ if $d \equiv 3$ (mod 4).   Define
\begin{displaymath}
  j(\gamma, z) := \epsilon_d^{-1} \left(\frac{c}{d}\right) \left(\frac{|cz+d|}{cz+d}\right)^{-1/2}
\end{displaymath}
for $\gamma = \left(\begin{smallmatrix}a & b\\ c & d\end{smallmatrix}\right) \in \Gamma_0(N)$ and $\Im z > 0$. Here $\left(\frac{c}{d}\right)$ is the extended Kronecker symbol as in \cite{Sh}. If $\chi$ is even, let $\textbf{H}_{\kappa}(N, \chi)$ be the Hilbert space of $L^2$-integrable functions $f$ satisfying
\begin{displaymath}
  f(\gamma z) = \chi(\gamma)  j(\gamma, z)^{2\kappa} f(z)
\end{displaymath}
for all $\gamma  \in \Gamma_0(N)$. For $t \in \Bbb{C}$, denote by $\textbf{H}_{\kappa}(N, \chi, t)$ the subspace of smooth functions $u \in \textbf{H}_{\kappa}(N, \chi)$ satisfying $(\Delta_{\kappa} + 1/4 +  t^2)u=0$. Without loss of generality we shall always assume $\Im t \geq 0$. 

For each equivalence class $\mathfrak{a}$ of cusps of $\Gamma_0(N)$ let $\Gamma_{\mathfrak{a}} := \{\gamma \in \Gamma_0(N) \mid \gamma\mathfrak{a} =\gamma\}$ be the stabilizer of $\mathfrak{a}$, $\sigma_{\mathfrak{a}} \in SL_2(\Bbb{R})$ be a scaling matrix (i.e. $\sigma_{\mathfrak{a}}\infty = \mathfrak{a}$ and $\sigma_{\mathfrak{a}}^{-1}\Gamma_{\mathfrak{a}} \sigma_{\mathfrak{a}} = \Gamma_{\infty}$) and $\gamma_{\mathfrak{a}} = \sigma_{\mathfrak{a}} \left(\begin{smallmatrix}1 & 1\\ & 1\end{smallmatrix}\right)\sigma_{\mathfrak{a}}^{-1} = \left(\begin{smallmatrix} \ast & \ast \\ c_{\mathfrak{a}} & d_{\mathfrak{a}} \end{smallmatrix} \right) \in \Gamma_0(N)$, say,  a generator of $\Gamma_{\mathfrak{a}}$. A cusp $\mathfrak{a}$ is singular for weight $\kappa$ and character $\chi$, if 
\begin{displaymath}
  \chi(d_{\mathfrak{a}}) \epsilon_{d_{\mathfrak{a}}}^{-2\kappa} \left(\frac{c_{\mathfrak{a}}}{d_{\mathfrak{a}}}\right) = 1.
\end{displaymath}
For a singular cusp let (initially for $\Re s > 1$)
\begin{displaymath}
  E_{\mathfrak{a}}(z; s) := \sum_{\gamma \in \Gamma_{\mathfrak{a}}\backslash \Gamma_0(N)} \bar{\chi}(\sigma_{\mathfrak{a}}^{-1}\gamma)  j(\sigma_{\mathfrak{a}}^{-1}\gamma, z)^{-2\kappa} \Im(\sigma_{\mathfrak{a}}^{-1}\gamma z)^s 
\end{displaymath}
be the Eisenstein series attached to $\mathfrak{a}$. We write the Fourier expansion as
\begin{displaymath}
  E_{\mathfrak{a}}(z; s)  = \delta_{\mathfrak{a} = \infty}y^s + \phi_{\mathfrak{a}}(0, s)y^{1-s} + \pi^s e\left(-\frac{\kappa}{4}\right)\sum_{n \not=0} |n|^{s-1}\phi_{\mathfrak{a}}(n, s) \frac{W_{\sgn(n)\frac{\kappa}{2}, \frac{1}{2}-s}(4\pi|n|y)}{\Gamma(s+\sgn(n)\frac{\kappa}{2})} e(nx)
\end{displaymath}
where $W_{\alpha, \beta}(y)$ is the standard Whittaker function, and 
\begin{displaymath}
  \phi_{\mathfrak{a}}(n, s) = \phi_{\mathfrak{a}}(n, s, \kappa, \chi) = \sum_{\substack{0 \leq d < c\\ \left(\begin{smallmatrix}* & *\\ c & d\end{smallmatrix}\right)   \sigma_{\mathfrak{a}}^{-1}\Gamma_0(N)}} \bar{\chi}(d) \left(\frac{c}{d}\right) \epsilon_d^{2\kappa} e\left(\frac{nd}{c}\right) c^{-2s}
\end{displaymath}
for $n \not=0$, cf.\ \cite[p.\ 3876]{Pr}. Note that
\begin{equation}
\phi_{\mathfrak{a}}(-n, s, \kappa, \chi) = \overline{\phi_{\mathfrak{a}}(n, \bar{s}, 2-\kappa, \bar{\chi})}.
\end{equation}
Let
\begin{displaymath}
 u_j(z) = \rho_j(0, y) + \sum_{n \not=0} \rho_j(n) W_{\sgn(n)\frac{\kappa}{2}, it_j}(4\pi |n|y)e(nx) \in \textbf{H}_{\kappa}(N, \chi)
\end{displaymath}
be a complete orthonormal set of automorphic eigenfunctions of $\Delta_{\kappa}$ with eigenvalues $\lambda_j = \frac{1}{4} + t_j^2$, i.e. $(\Delta_{\kappa} + \lambda_j)u_j = 0$.  We call $\lambda_j$ exceptional if $t_j \not\in \Bbb{R}$, i.e.\ $\lambda_j < 1/4$. The functions $u_j$  may be cusp forms (in which case $\rho_j(0, y) = 0$) or residues of possible poles of an Eisenstein series $E_{\mathfrak{a}}(z; s)$.  In either case, $\lambda_j \geq 3/16$ (cf.\ the discussion preceding (2.8) below). 
The space $\textbf{H}_{\kappa}(N, \chi, i/4)$ corresponding to the exceptional eigenvalue $3/16$ is the kernel of $\Lambda_{\kappa}$. Hence if $(\Delta_{\kappa} + 3/16)u = 0$, then $y^{-\kappa/2}$ is holomorphic, and so $\textbf{H}_{1/2}(N, \chi, i/4) = \{y^{1/4}f \mid f \in M_{1/2}(N, \chi)\}$ and $\textbf{H}_{3/2}(N, \chi, i/4) = \{y^{3/4}f \mid f \in S_{3/2}(N, \chi)\}$ where $M_{\kappa}(N, \chi)$ and $S_{\kappa}(N, \chi)$ are the spaces of holomorphic modular forms of weight $\kappa$, level $N$ and character $\chi$. 
In particular,  $u \in H_{\kappa}(N, \chi, i/4)$ has no  negative Fourier coefficients.

If $u$ is any automorphic eigenfunction of $\Delta_{\kappa}$ with spectral parameter $t = \sqrt{\lambda-1/4}$, then its Shimura lift is an even weight Maa{\ss} form with spectral parameter $2t$ (see e.g.\ \cite{Bi}). It is a cusp form  unless $u$ comes from theta-functions, so that $\lambda = 3/16$. In all other cases  the Kim-Sarnak bound \cite{KS} implies
\begin{equation}
   |\Im t | \leq \frac{\theta}{2} \leq \frac{7}{128}.
\end{equation}   
Let $t \not=i/4$.   As is \cite[p.\ 507 and 509]{DFI} we see that
\begin{equation}
 T_{\kappa} :=  \left(\frac{1}{16}+t^2\right)^{-1/2}X\Lambda_{\kappa} : \textbf{H}_{\kappa}(N, \chi, t) \rightarrow \textbf{H}_{2-\kappa}(N, \chi, t)
\end{equation}
is a bijective isometry. Essentially, $T_{\kappa}$ interchanges positive and negative Fourier coefficients; more precisely, combining \cite[(4.17), (4.18), (4.27), (4.28), 4.64)]{DFI} (which also hold for non-integral weight), we see that
\begin{equation}
\begin{split}
  T_{\kappa} & \left(\sum_{n \not=0} \rho(n) e(nx)W_{\sgn(n)\frac{\kappa}{2}, it}(4 \pi|n|y) \right) \\
  & = \sum_{n \not= 0} \rho(-n) \sgn(n)  \left(\frac{1}{16}+t^2\right)^{-\sgn(n)/2} e(nx)W_{\sgn(n)\frac{2-\kappa}{2}, it}(4 \pi|n|y). 
  \end{split}
\end{equation}

For $5/2 \leq \ell \in \frac{1}{2}\Bbb{Z} \backslash\Bbb{Z}$ let $S_{\ell}(N, \chi)$ denote the (finite-dimensional) Hilbert space of holomorphic cusp forms of weight $\ell$ and character $\chi$ for $\Gamma_0(N)$. For each $\kappa $ we choose an $L^2$-orthonormal basis
 \begin{displaymath}
  f_{\ell, j}(z) = \sum_{n > 0} (4\pi n)^{\ell/2} \rho_{\ell, j}(n) e(nz) \in S_{\ell}(N, \chi), \quad 1 \leq j \leq D_{\ell} := \text{dim}_{\Bbb{C}}S_{\ell}(N, \chi).
\end{displaymath}
Let $\phi$ be a smooth function on $[0, \infty)$ such that $\phi(0) = \phi'(0) = \phi''(0)=0$ and $\phi^{(j)}(x) \ll x^{-2-\varepsilon}$ for $0 \leq j \leq 3$ and $x \rightarrow \infty$. Assume that $\phi' \in L^{1}(\Bbb{R}_{>0}, dx/x)$. Define
\begin{displaymath}
\begin{split}
  \tilde{\phi}(t) :=  &\int_0^{\infty} J_{t-1}(x)\phi(x) \frac{dx}{x},\\
   \widehat{\phi}(t) := &\frac{\pi^2 e((\kappa+1)/4)}{\sinh(\pi t)(\cosh(2\pi t)+ \cos(\pi\kappa))\Gamma(\frac{1-\kappa}{2}+it)\Gamma(\frac{1-\kappa}{2}-it)}\\ & \times  \int_0^{\infty}\left( \cos\left(\pi\left(\kappa/2+it\right)\right) J_{2it}(x)-\cos\left(\pi\left(\kappa/2-it\right)\right) J_{-2it}(x)\right)\phi(x)\frac{dx}{x}\\
   \check{\phi}(t) := & 2 e\left(\frac{\kappa}{4}\right) \cosh(\pi t) \int_0^{\infty} K_{2it}(y) \phi(y) \frac{dy}{y}.
\end{split}  
\end{displaymath}
Here $K_{\nu}(x)$ is the Bessel $K$-function given by
\begin{equation}
\begin{split}
   K_{\nu}(x) & = \frac{\pi}{\sin(\pi \nu)} (I_{-\nu}(x)-I_{\nu}(x)) \\
   & = \frac{\pi}{\sin(\pi \nu)} \sum_{j=0}^{\infty}\left(\frac{(x/2)^{-\nu+2j}}{j!\Gamma(j+1-\nu)}- \frac{(x/2)^{\nu+2j}}{j!\Gamma(j+1+\nu)} \right), \quad x \not\in (-\infty, 0].
   \end{split}
\end{equation}
We observe that
\begin{equation}
  \widehat{\phi}(i/4) = \left\{\begin{array}{ll}
    (1+i)\int_{0}^{\infty} \cos(x) \phi(x) x^{-3/2} dx,  & \kappa = 1/2,\\
     (i-1)/2 \int_{0}^{\infty} \sin(x) \phi(x) x^{-3/2} dx, & \kappa = 3/2. \end{array} \right.
\end{equation}
Finally,  if $4 \mid N \mid c$, let 
\begin{equation}
  K_{\chi, \kappa}(m, n; c) :=  \left\{ \begin{array}{ll}
    \sum_{d \, (c)}^{\ast}  \chi(d) \epsilon_d \left(\frac{c}{d}\right) e\left(\frac{m d + n\bar{d}}{c}\right), & \kappa = 1/2\\
    \sum_{d \, (c)}^{\ast}  \chi(d) \epsilon_d^{-1} \left(\frac{c}{d}\right) e\left(\frac{m d + n\bar{d}}{c}\right), & \kappa = 3/2
    \end{array}\right.
\end{equation}
denote the (twisted)  Kloosterman sum with theta-multiplier. Note that
\begin{equation}
  K_{\chi, \kappa}(-m, -n;c) = \overline{K_{\bar{\chi}, 2-\kappa}(m, n; c)}
\end{equation}
 It satisfies a Weil-type bound (see \cite[(1.6)]{Iw}, or \cite[ch.\ 12.3]{IK} for the underlying theory)
\begin{equation}
  |K_{\chi, \kappa}(m, n; c)| \leq \tau(c) (m, n, c)^{1/2} c^{1/2}. 
\end{equation}
Note that the definition (2.13) makes sense for odd characters as well, and we have
\begin{equation}
  K_{\chi, \kappa}(m, n;c) = K_{\chi \chi_{-4}, 2-\kappa}(m, n; c)
\end{equation}
where $\chi_{-4} = \left(\frac{-4}{.}\right)$. With these preparations we can state the Kuznetsov sum formula as generalized by Proskurin for half-integral weight. 
\begin{prop} [Kuznetsov-Proskurin] For $m, n \geq 1$ and $\chi$ an even character one has, with the above notations and assumptions,
\begin{displaymath}
\begin{split}
 & \sum_{N \mid c} \frac{1}{c} K_{\chi, \kappa}(m, n; c) \phi\left(\frac{4\pi \sqrt{mn}}{c}\right) \\
 & = 4\sqrt{mn}\left(\left.\sum_j\right.^{(\kappa)}  \frac{\bar{\rho}_j(m)\rho_j(n) }{\cosh(\pi t_j)} \widehat{\phi}(t_j) + \sum_{\substack{\ell \geq 5/2\\ \ell \equiv \kappa \, (\text{mod }2)}}  \Gamma(\ell)e(\ell/4)\tilde{\phi}(\ell)\sum_{j=1}^{D_{\ell}} \bar{\rho}_{\ell, j}(m)\rho_{\ell, j}(n)\right) \\
  & + \sum_{\mathfrak{a} \text{ singular}} \int_{\Bbb{R}} \left(\frac{n}{m}\right)^{it} \frac{\bar{\phi}_{\mathfrak{a}}(m, 1/2 + it)\phi_{\mathfrak{a}}(n, 1/2+it)}{\cosh(\pi t)|\Gamma(\frac{\kappa+1}{2}+ it)|^2} \widehat{\phi}(t)dt 
     \end{split}
\end{displaymath}
as well as
\begin{displaymath}
\begin{split}
 & \sum_{N \mid c} \frac{1}{c} K_{\chi, \kappa}(m, -n; c) \phi\left(\frac{4\pi \sqrt{mn}}{c}\right)  = 4\sqrt{mn}\left.\sum_j\right.^{(\kappa)}  \frac{\bar{\rho}_j(m)\rho_j(-n) }{\cosh(\pi t_j)} \check{\phi}(t_j)  \\
& \quad\quad + \sum_{\mathfrak{a} \text{ singular}} \int_{\Bbb{R}} \left(\frac{n}{m}\right)^{it} \frac{\bar{\phi}_{\mathfrak{a}}(m, 1/2 + it)\phi_{\mathfrak{a}}(-n, 1/2+it)}{\cosh(\pi t)\Gamma(\frac{\kappa+1}{2}- it)\Gamma(\frac{1-\kappa}{2}+it)} \check{\phi}(t)dt. 
     \end{split}
\end{displaymath}
Here $\sum^{(\kappa)}$ indicates that the sum is over an orthonormal system in $\textbf{H}_{\kappa}(N, \chi)$.
\end{prop}

\textbf{Proof.} The first formula is \cite[Theorem]{Pr}.  The second formula does not seem to be in the literature for non-zero weight, but can be proved along the same lines (cf.\ \cite{Pr, DI, Mo}). We postpone the proof to section 4. \\

We shall need the first formula of the previous proposition in a somewhat less refined, but more explicit version \cite[Lemma 3]{Pr}.

\begin{lemma} Let $r \in \Bbb{R}$, and define
\begin{displaymath}
  H(t, r) := \frac{\cosh(\pi t)\Gamma(\frac{1}{2} + i(r-t))\Gamma(\frac{1}{2} + i(r+t))\Gamma(\frac{1}{2} - i(r-t))\Gamma(\frac{1}{2} - i(r+t)) }{2\Gamma(1-\frac{\kappa}{2}+it)\Gamma(1-\frac{\kappa}{2}+it)}. 
\end{displaymath}
For $m, n \geq 1$ one has
\begin{displaymath}
\begin{split}
 & \sum_{N \mid c} \frac{4\sqrt{mn}}{c^2} K_{\chi, \kappa}(m, n; c) e\left(-\frac{\kappa+1}{4}\right)\int_{-i}^i K_{ir}\left(\frac{4\pi\sqrt{mn}}{c}y\right)y^{\kappa-1}dy + \frac{\delta_{mn}}{2\pi} \\
 & = 4\sqrt{mn}\left.\sum_j\right.^{(\kappa)}  \frac{\bar{\rho}_j(m)\rho_j(n) }{\cosh(\pi t_j)}  H(t_j, r)+ \sum_{\mathfrak{a}\text{ singular}} \int_{\Bbb{R}} \left(\frac{n}{m}\right)^{it} \frac{\bar{\phi}_{\mathfrak{a}}(m, 1/2 + it)\phi_{\mathfrak{a}}(n, 1/2+it)}{\cosh(\pi t)|\Gamma(\frac{\kappa+1}{2}+ it)|^2} H(t, r)dt
       \end{split}
       \end{displaymath}
 where the integration on the left hand side is counterclockwise over the right half of the unit circle.      
\end{lemma}
It should be noted that in our case we have Weil's bound (2.15) for the Kloosterman sums in question, so the left hand side converges absolutely. More precisely, Lemma 1 and 2 in \cite{Pr} hold in $\Re s > 3/4$, and so we can apply Lemma 3 in \cite{Pr} with $\sigma = 1$. 

The holomorphic analogue of the preceding formulas reads as follows (see \cite[p.\ 389]{IK}).

\begin{lemma} [Petersson] For $m, n \geq 1$, $5/2 \leq \ell \in \frac{1}{2}\Bbb{Z}\setminus\Bbb{Z}$, $\kappa \equiv \ell$ (mod 2),  one has
\begin{displaymath}
4\sqrt{mn}\Gamma(\ell)e(\ell/4)\sum_{j=1}^{D_{\ell}} \bar{\rho}_{\ell, j}(m)\rho_{\ell, j}(n)=  (\ell-1)\left(\frac{\delta_{mn} e(\ell/4)}{\pi} + 2 \sum_{N \mid c}\frac{1}{c}  K_{\chi, \kappa}(m, n; c)J_{\ell-1}\left(\frac{4 \pi \sqrt{mn}}{c}\right)\right).
\end{displaymath}
\end{lemma}

Lemmas 3 and 4 together with (2.15) give bounds for Fourier coefficients:

\begin{lemma} Let  $n$ be a positive integer, $T \geq 1$.   Then
\begin{displaymath}
\begin{split}
 \left. \sum_{|t_j| \leq T}\right.^{(\kappa)} \frac{n |\rho_j(\pm n)|^2(1+|t_j|)^{\pm \kappa-1/2}}{\cosh(\pi t_j)} \\
 + \sum_{\mathfrak{a}\text{ singular}} \int_{-T}^T \frac{ |\phi_{\mathfrak{a}}(\pm n, 1/2 + it)|^2(1+|t|)^{\pm \kappa-1/2}}{\cosh(\pi t)|\Gamma(\pm \frac{k}{2} + \frac{1}{2}+it)|^2} dt &\ll T^{3/2} + (nT)^{1/2+\varepsilon};\\ 
   \sum_{\substack{\ell \leq 5/2\\ \ell \equiv  \kappa \, (2)}}  \sum_{j=1}^{D_{\ell}} n |\rho_{\ell, j}(n)|^2& \ll n^{1/2+\varepsilon}
  \end{split}
\end{displaymath}
\end{lemma}

\textbf{Proof.} For the first bound with the $+$ sign, we multiply the formula in Lemma 3 (with $m=n$) by 
\begin{displaymath}
  (1+|r|)^{1/2} e^{-(r/T)^2} 
\end{displaymath}  
and integrate over $r$. By Stirling's formula,
\begin{displaymath}
  \int_{\Bbb{R}} H(t, r) (1+|r|)^{1/2} e^{-(r/T)^2} dr \gg \int_{|t|}^{|t|+1} r^{\kappa-1/2} e^{-(r/T)^2} dr \gg (1+|t|)^{\kappa-1/2} e^{-(|t|/T)^2}, 
\end{displaymath}
so the right hand side becomes an upper bound for the quantity we want to estimate.  Now we estimate the integral over the left side in Lemma 3. The diagonal term contributes $O(T^{3/2} )$. For the other term we substitute the integral representation \cite[8.432.1]{GR}
\begin{displaymath}
  K_{2 ir}(y)  = \int_0^{\infty} e^{-y \cosh x} \cos(2 r x) dx,
\end{displaymath}
so that by (2.15) we need to estimate
\begin{displaymath}
\begin{split}
&  n \sum_c \frac{(c, n)^{1/2}}{c^{3/2-\varepsilon}}   \int_{-i}^i \int_{-\infty}^{\infty}\int_0^{\infty} (1+|r|)^{1/2}e^{-(r/T)^2} \exp\left(-\frac{4\pi y n\cosh(x)}{c}\right)\cos(2rx)\, dx\, dr \, y^{\kappa-1}dy.
 \end{split}
\end{displaymath}
The integral over $r$ is, by repeated partial integration, at most $\ll_A T^{3/2} (Tx)^{-A}$ for any $A \geq 0$. Now the $x$-integral is at most $\ll e^{- 4 \pi n (\Re y)/c} T^{1/2}$, and so the $y$-integral is $\ll T^{1/2} \min(1, c/n) \ll T^{1/2} (c/n)^{1/2-2\varepsilon}$.  Summing over $c$, we  obtain a total contribution of $\ll (Tn)^{1/2+\varepsilon}$. Observing (2.7) and Stirling's formula, we can also include negative Fourier coefficients of Eisenstein series into the result. In order to include negative cusp form coefficients, we first observe that the spectral parameter $t=i/4$ does not contribute to this sum, since all functions in $\textbf{H}_{\kappa}(N, \chi, i/4)$ have no negative Fourier coefficients. Using (2.9) and (2.10) we see that
\begin{displaymath}
\begin{split}
  &\left. \sum_{|t_j| \leq T}\right.^{(\kappa)} \frac{n |\rho_j(-n)|^2(1+|t_j|)^{- \kappa-1/2}}{\cosh(\pi t_j)} = \left. \sum_{|t_j| \leq T}\right.^{(2-\kappa)} \frac{n |\rho_j(n)|^2 (\frac{1}{16}+t_j^2)(1+|t_j|)^{- \kappa-1/2}}{\cosh(\pi t_j)}\\
  & \ll \left. \sum_{|t_j| \leq T}\right.^{(2-\kappa)} \frac{n |\rho_j(n)|^2 (1+|t_j|)^{2- \kappa-1/2}}{\cosh(\pi t_j)},
\end{split}
\end{displaymath}
and for the latter sum we have already proved the desired bound.

For the second estimate, we use Lemma 4 for each $\ell$, so the term in question is at most
\begin{displaymath}
 \ll \sum_{\ell \geq 5/2} \frac{\ell}{\Gamma(\ell)}\left(1 + \sum_{N \mid c}\frac{(n, c)\tau(c)}{\sqrt{c}} J_{\ell-1}\left(\frac{n}{c}\right)  \right)\ll n^{1/2+\varepsilon}
\end{displaymath}
using the uniform bound $J_{\ell - 1}(x) \leq 1$ which follows, for example, from the integral representation \cite[8.411.1]{GR}. \\



Next we collect various bounds for the integral transforms appearing in the Kuznetsov formula.
\begin{lemma} Let $Y \leq 1 \leq P$, and let $\phi$ be a function supported on $[Y/2, 2Y]$ such that $\phi^{(j)} \ll_j (P/Y)^j$ for all $j \in \Bbb{N}_0$. Then
\begin{displaymath}
\begin{split}
\tilde{\phi}(t) & \ll \Gamma(t)^{-1}, \quad t \geq 1,\\
\widehat{\phi}(t) & \ll_A  Y^{-2 |\Im t|} (1+ |t|)^{\kappa-\frac{1}{2}}\left(1+\frac{|t|}{P}\right)^{-A}, \quad t \in [0, \infty) \cup [0, i/4],\\
\check{\phi}(t) & \ll_A Y^{- 2 |\Im t|} (1+\abs{t})^{-1/2}\left(1+\frac{|t|}{P}\right)^{-A} , \quad t \in [0, \infty) \cup [0, i/4], \\
\widehat{\phi}(i/4) & \ll Y^{1/2}, \quad \kappa=3/2
\end{split}
\end{displaymath}
for any $A \geq 0$. Moreover,  if $\kappa=1/2$ and   $\phi$ is any function (satisfying the hypotheses of Proposition 2), then
\begin{displaymath}
 \widehat{\phi}(i/4) = (1+i)\int_{0}^{\infty} \phi(x) x^{-3/2} dx +O\left( \int_{0}^{\infty} \min(1, x^2) |\phi(x)| x^{-3/2} dx\right). 
\end{displaymath}

\end{lemma}

\textbf{Proof.} This follows easily from  (2.2), (2.11) and repeated integration by parts, and by (2.12) noting that  $|\sin(x)| \leq x$ and  $\cos(x) = 1+O( \min(1, x^2))$. \\

 

Finally we evaluate the crucial character sum.

\begin{lemma} Let $(d, c) = 1$, $u \in \Bbb{Z}$, $r \in \Bbb{N}$, and assume $4 \mid c$. Then
\begin{displaymath}
  G(d, u; c) := \sum_{b \, (c)} e\left(\frac{db^2 + ub}{c}\right) = \left\{\begin{array}{ll}
  0, & 2 \nmid u\\
  (1+i)  \sqrt{c} \left(\frac{c}{d}\right)\epsilon_d^{-1} e\left(\frac{-\bar{d} u^2/4}{c}\right), & 2 \mid u.
  \end{array}\right.
\end{displaymath} 
\end{lemma}

\textbf{Proof.} This is a special case of e.g.\  \cite[Lemma 2]{Bl}.

\section{Proof of Theorem  2}

Let $X, P \geq 1$, $f \in S_k(N, \chi)$, $q(x) = x^2 + sx+t \in \Bbb{Z}[x]$ with discriminant $\Delta = s^2-4t$ and the smoothing function $w$ be as in Theorem 2. Of course, $s$ and $t$ are unrelated to the variable $s$ and the spectral parameter $t$ in the previous section. In the following all implied constants may depend on $f$ and $q$. We start with a few preliminary reductions.  

Let us first assume that $\Delta=0$. Then we are essentially in a symmetric square situation. More precisely, $s = 2s'$ is necessarily even, and $q(x) = (x+s')^2$.  We can assume that the cusp form $f$ in Theorem 2 is of the form $f(z) = g(dz)$ for some newform $g$ of level dividing $N/d$.  Let $\lambda_g(n)$ denote the Hecke eigenvalues of $g$ and write $d = d_1d_2^2$ with $\mu^2(d_1)=1$. Then
\begin{displaymath}
  \sum_n a(q(n)) w(n) = \sum_m \lambda_g(d_1m^2) w(d_1d_2m - s').
\end{displaymath}
For simplicity, let us write $v(m) := w(d_1d_2m - s')$, and denote by $\hat{v}(s) = \int_0^{\infty} v(y) y^{s-1}dy$ the Mellin transform of $v$. By partial integration, $\hat{v}(s) \ll X^{\Re s} (P/(|s|+1))^{2}$ in fixed vertical strips. Now the Hecke relation $\lambda(ab) = \sum_{d \mid (a, b)} \mu(d)\chi(d) \lambda(a/d)\lambda(b/d)$ yields
\begin{displaymath}
  \sum_m \frac{\lambda_g(d_1m^2)}{m^s} =\sum_{f \mid d_1} \frac{\mu(f)\chi(f)\lambda_g(d_1/f)}{f^s}  \sum_m \frac{\lambda_g(fm^2)}{m^s},
\end{displaymath}
from which we readily obtain  
\begin{displaymath}
  \sum_m \frac{\lambda_g(d_1m^2)}{m^s}=\lambda_g(d_1) \prod_{p \mid d_1} \left(1+\frac{\chi(p)}{p^s}\right)^{-1} \sum_m \frac{\lambda_g(m^2)}{m^s} .
\end{displaymath}
Thus
 \begin{displaymath}
 \begin{split}
 \sum_n a(q(n)) w(n)  &=    \frac{1}{2\pi i} \int_{(2)} \sum_m \frac{\lambda_g(d_1m^2)}{m^s} \hat{v}(s) ds \\
 & = \frac{1}{2\pi i} \int_{(2)} \lambda_g(d_1) \prod_{p \mid d_1} \left(1+\frac{\chi(p)}{p^s}\right)^{-1} L(s, \text{sym}^2 g) L(2s, \chi^2)^{-1} \hat{v}(s) ds.
 \end{split}
\end{displaymath}
It is known \cite[p.\ 95, Remark 2]{Sh2} that $L(s, \text{sym}^2g)$ is holomorphic in $\Re s \geq 1/2$ except for a possible pole at $s=1$ which can only occur if $\chi$ is quadratic and $k$ is odd. Shifting the line of integration to $\Re s = 1/2$ and using the convexity bound $L(\text{sym}^2g, 1/2 + it) \ll (|t|+1)^{3/4+\varepsilon}$ we get the bound
\begin{equation}
\begin{split}
  & \sum_n a(q(n)) w(n)\\
  & = \frac{\lambda_g(d_1)}{d_1d_2} \prod_{p \mid d_1} \left(1+\frac{\chi(p)}{p}\right)^{-1}   \frac{\text{res}_{s=1} L(s, \text{sym}^2 g)}{L(2, \chi^2)}\int_{0}^{\infty} w(y) dy  + O( X^{1/2} P^{2})
\end{split}
\end{equation}
which yields Theorem 2 in the case $\Delta = 0$. 

From now on we assume $\Delta \not= 0$.  Let $N' := [4, N]$. We view $f \in S_k(N, \chi) \subseteq S_k(N', \chi)$ as a form of level $N'$. 
Assume first\footnote{Such an assumption can certainly not be satisfied for newforms.} that $a(n) = 0$ for all $n$ unless $n =  \square$. Then $a(q(n)) = 0$, unless $(2n+s)^2 - \Delta = \square$. For $\Delta \not=0$ there are only finitely many $n$ of this kind, so we are done in this case. Let $S^{\ast}_k(N', \chi) = \langle P_m \mid m \not= \square \rangle$. We decompose $f=f^{\ast} + f^{\perp}$ where $f^{\ast} \in S^{\ast}_k(N', \chi)$, and $f^{\perp}$ is in the orthogonal complement. Then by (2.3) all coefficients of $f^{\perp}$ vanish except those with square index, so its contribution is negligible. If $f$ does not have the property that $a(n) = 0$ unless $n =  \square$, then $f$ is in the space generated by all $P_m$ with $m \not= \square$: indeed, if  $f$ was orthogonal to this space, then by (2.2) all coefficients $a(n)$ vanish except if $n =  \square$. 
Therefore we can assume that $f$ is a finite linear combination of certain $P_m$ with $m\not= \square$ getting 
\begin{displaymath}
  \sum_{n} a(q(n)) w(n) = \sum_{\substack{m \ll 1\\m \not= \square}} \alpha_m \sum_n w(n) \sum_{N' \mid c} \frac{1}{c} S_{\chi}(m, q(n); c) J_{k-1}\left(4 \pi\frac{\sqrt{q(n)m}}{c}\right) +O(1)
\end{displaymath}
for some $\alpha_m \in \Bbb{C}$. We could have avoided this little maneuver of excluding square $m$, but it simplifies later calculations a bit. Opening the Kloosterman sum, the absolutely convergent double sum over $n$ and $c$ now equals
\begin{displaymath}
   \sum_{N' \mid c} \frac{1}{c} \left. \sum_{d \, (c)}\right. ^{\ast} \chi(d) e\left(\frac{md}{c}\right) \sum_n w(n) e\left(\frac{q(n)\bar{d}}{c}\right) J_{k-1}\left(\frac{4\pi\sqrt{q(n)m}}{c}\right).
\end{displaymath}
We split the inner sum into residue classes modulo $c$ and apply Poisson summation (Lemma 1) getting
\begin{equation}
 \sum_h  \sum_{N' \mid c} \frac{1}{c^2} \left. \sum_{d \, (c)}\right. ^{\ast} \chi(d) e\left(\frac{md}{c}\right) \sum_{b \, (c)} e\left(\frac{q(b)\bar{d} + hb}{c}\right) g(h; c)
 \end{equation}
where
\begin{displaymath}
  g(h; c) := \int_{\Bbb{R}} w(y) J_{k-1}\left(\frac{4\pi\sqrt{q(y)m}}{c}\right) e\left(-\frac{hy}{c}\right) dy.
\end{displaymath}
Integrating by parts twice and using (2.4), (2.5) and Lemma 7, one easily sees that for $k \geq  4$ the double sum over $h$ and $c$ is absolutely convergent, thereby justifying the change of summation. For later purposes, let us define
\begin{equation}
    g^{\ast}(h; c) := \delta(h) \int_{\Bbb{R}} w\left(y-\frac{s}{2}\right) J_{k-1}\left(\frac{4\pi y \sqrt{m}}{c}\right) \cos\left(\frac{2\pi hy}{c}\right) dy, \quad c >0, h \in \Bbb{Z}, 
\end{equation}
where $\delta(h) = 2$ if $h \not= 0$ and $\delta(0) = 1$. Then (recall $q(x) = x^2 + sx + t$ and $\Delta = s^2-4t$)
\begin{displaymath}
\begin{split}
 \sum_{\pm} & g(\pm h; c)e\left(\mp \frac{sh}{2}\right)  =\delta(h) \int_{\Bbb{R}} w\left(y-\frac{s}{2}\right) J_{k-1}\left(\frac{4\pi \sqrt{m(y^2-\Delta/4)}}{c}\right)\cos\left(\frac{2\pi hy}{c}\right)dy\\
  &  = g^{\ast}(h; c) +   O\left( \int_{\Bbb{R}} w\left(y-\frac{s}{2}\right) \left(J_{k-1}\left(\frac{4\pi \sqrt{m(y^2-\frac{\Delta}{4})}}{c}\right)-J_{k-1}\left(\frac{4\pi \sqrt{my^2}}{c}\right)\right) \cos\left(\frac{2\pi hy}{c}\right)dy\right)
\end{split}
\end{displaymath}
where $\sum_{\pm} f(\pm 0)$ is interpreted as $f(0)$. Integrating by parts twice and using the mean value theorem, (2.4) and (2.5) we find
\begin{equation}
  \sum_{\pm}  g(\pm h; c)e\left(\mp \frac{sh}{2}\right) =  g^{\ast}(h; c) +   O\left( \frac{1}{h^2c^{1/10}} \left(\frac{P^2}{X^{9/10}} + 1\right)\right)
\end{equation}
for $k \geq 4$. 

\begin{lemma} a)  Let $0 < \varepsilon < 1/50$. Then $g^{\ast}(h;c) \ll_{\varepsilon} \min(1, \sqrt{c}) (1+|h|)^{-5/2}X^{-10}$ unless
\begin{displaymath}
  \frac{(1+|h|)}{c} \leq  X^{\varepsilon} \frac{P}{X}.
 \end{displaymath}
b) One has $g^{\ast}(h; c) \ll X (X/c)^3 \ll  X^{-26}(1+|h|)^{-3}$ unless $c \leq X^{10} (1+|h|)$.  
 \end{lemma}

\textbf{Proof.} a) Let us first assume  $c \leq X^{1-\varepsilon/2}P^{-1}$. If $|h|$ is very large, say $|h| \geq X^{100}$, we integrate by parts in (3.3) three times the cosine factor and differentiate the Bessel function using (2.5) and (2.4). In this way we get the bound
\begin{displaymath}
  X \left(\frac{X}{c}\right)^{-1/2}\left(\frac{P}{X} + \frac{1}{c}\right)^3 \left(\frac{c}{|h|}\right)^3
\end{displaymath}
which is acceptable.  Otherwise we can insert the asymptotic  expansion (2.6) into (3.3) getting an oscillating factor 
\begin{displaymath}
  \phi(y) = e\left(\frac{  y (\pm 2 \sqrt{m} \pm h)}{c}\right)
\end{displaymath}
for various choices of $\pm$ in the integral. Choosing $A$ in (2.6) large enough, the error is admissible. Now we integrate by parts $\lceil \frac{22}{\varepsilon} + 4\rceil$ times. Since $m \not= \square$ and $m$ is bounded, each integration of $\phi$ introduces at most a factor $\ll c/(1+|h|)$, while each differentiation of $w$ introduces a factor $P/X$. Therefore $g^{\ast}(h; c) \ll_{\varepsilon} \min(1, c) (1+|h|)^{-5/2} X^{-10}$. 

Let us now assume $c \geq X^{1-\varepsilon/2}P^{-1}$. We have also the general assumption $(1+|h|)/c > PX^{\varepsilon-1}$, for otherwise the statement of the lemma is void. These two bounds imply $(1+|h|) \geq X^{\varepsilon/2}$.  Now we integrate by parts in (3.3), this time integrating only the cosine factor. Hence we can bound (3.3) by
\begin{displaymath}
  X\left(\frac{c}{1+|h|}\left(\frac{P}{X}+ \frac{1}{c}\right)\right)^{\lceil 162/\varepsilon \rceil} \ll \frac{X}{X^{162}} + \frac{X}{(1+|h|)^{162/\varepsilon}} \ll 
  \frac{1}{(1+|h|)^{5/2}X^{10}} 
\end{displaymath}
if  $1+|h| \leq X^{60}$. If $1+|h| > X^{60}$ and $c \leq (1+|h|)^{11/12}$, we estimate the left hand side of the preceding display by
\begin{displaymath}
    X\left(\frac{c}{1+|h|}\left(\frac{P}{X}+ \frac{1}{c}\right)\right)^{\lceil 162/\varepsilon \rceil} \ll \frac{X}{(1+|h|)^{13/\varepsilon}  }
\end{displaymath} 
which is also acceptable. 
Finally, if $1+|h| > X^{56}$ and $c > (1+|h|)^{11/12}$, we estimate (3.3) trivially by (2.4) getting the bound
\begin{displaymath}
  \ll X \left(\frac{X}{c}\right)^{k-1} \leq \frac{X^4}{(1+|h|)^{11/4}} \leq \frac{1}{(1+|h|)^{5/2}X^{10}}
\end{displaymath}
for $k \geq 4$. 

b) This follows directly from (2.4) if $k \geq 4$. This completes the proof of the lemma. \\

\textbf{Remark.} Of course, the exponents in the statement and the proof of the lemma are fairly arbitrary; the only constraint in part a) is that our method would give at most a saving of $1/(1+|h|)^{k-1-\varepsilon} \ll 1/(1+|h|)^{3-\varepsilon}$  for $k \geq 4$ in the $h$-aspect. \\

We now evaluate the $b$-sum in (3.2) by Lemma 7. Using the definition (2.13), we can recast (3.2) as
\begin{equation}
\begin{split}
&  (1+i)  \sum_{h \equiv s \, (2)}  \sum_{N' \mid c} \frac{g(h;c)}{c^{3/2}}  \left. \sum_{d \, (c)}\right. ^{\ast}  \chi(d) e\left(\frac{md+t\bar{d}}{c}\right)  \left(\frac{c}{d}\right)\epsilon_d^{-1}e\left(\frac{-d(s\bar{d}+h)^2/4}{c}\right)\\
=&  \frac{(1+i)}{4} \sum_{h \equiv s \, (2)}  \sum_{4N' \mid c} \frac{g(h;c/4)}{(c/4)^{3/2}} e\left(-\frac{2sh}{c}\right)  K_{\chi, \frac{3}{2}}(4m-h^2, -\Delta, c)\\
= &  2(1+i) \sum_{\substack{h \equiv s \, (2)\\h \geq 0}}  \sum_{4N' \mid c} \frac{g^{\ast}(h;c/4)}{c^{3/2}}  K_{\chi, \frac{3}{2}}(4m-h^2, -\Delta, c) + O\left(\frac{P^{2}}{X^{9/10}} + 1\right)
\end{split}
\end{equation}
where we used (2.15) and (3.4) in the final step. By our preliminary remarks we can assume that $\Delta \not=0$ and $m \not=\square$, so neither of the entries of the Kloosterman sum vanishes. Let $\omega$ be a nonnegative smooth function such that $\omega =1$ on $[0, 1]$ and $\omega = 0$ on $[2, \infty)$. By (2.15), (3.3) and Lemma 8,  we can replace $g^{\ast}$ by   
\begin{displaymath}
  \tilde{g}(h; c) := g^{\ast}(h; c) \omega\left(\frac{(1+|h|)X}{cX^{\varepsilon}P}\right)\omega\left(\frac{c}{X^{10}(1+|h|)}\right)
\end{displaymath}
at the cost of an error of at most
\begin{displaymath}
\sum_h  \sum_{c \leq (1+|h|)X/P} \frac{\tau(c)}{c(1+|h|)^{5/2} X^{10}} + \sum_h \sum_{c \geq (1+|h|)X^{10}}  X\left(\frac{X}{c}\right)^{3} \ll \frac{1}{X^{9}}.
\end{displaymath}
In order to apply the Kuznetsov formula let us finally define
\begin{equation}
  \phi_h(z) :=  \frac{2(i+1)z^{1/2}}{|\Delta(4m-h^2)|^{1/4}} \tilde{g}\left(h; \frac{\sqrt{|\Delta(4m-h^2)|}}{4z}\right),
\end{equation}
then the main term on the right hand side of (3.5) is just
\begin{equation}
\sum_{\substack{h \equiv s \, (2)\\h \geq 0}}  \sum_{4N' \mid c} \frac{1}{c} K_{\chi, \frac{3}{2}}(4m-h^2, -\Delta, c) \phi_h\left(\frac{\sqrt{|\Delta(4m-h^2)|}}{c}\right),
\end{equation} 
cf.\ (1.6) in the introduction. We collect some properties of $\phi_h$.

\begin{lemma}
One has
\begin{displaymath}
\begin{split}
 &  \phi_h(z) = 0 \text{ if } z \gg P/X^{1-\varepsilon} \text{ or } z \ll X^{-10},\\
 & \phi_h^{(j)}(z) \ll_j \frac{z^{1/2}X}{(1+|h|)^{1/2}} \left(\frac{1+zX}{z}\right)^j \min\left( \left(\frac{zX}{1+|h|}\right)^{k-1}, \left(\frac{zX}{1+|h|}\right)^{-1/2}\right) =: \Xi(h) \left(\frac{1+zX}{z}\right)^j, \\
 & \int_0^{\infty} \phi_h(z) z^{-3/2} dz \ll X^{-9}(1+|h|)^{-9/4}, \text{ if } |h| > 2\sqrt{m} \text{ and } k \text{ is even.}
\end{split}
\end{displaymath}
\end{lemma}

\textbf{Remark.} An inspection of the proof of the last statement shows that the condition $|h| > 2 \sqrt{m}$ comes from an identity among special functions that in our context captures aritmetic information. It will turn out to be significant for the rest of the proof. \\

\textbf{Proof.} The first statement follows directly from the properties of $\omega$. To see the second statement, use (2.4), (2.5) and distiguish the cases $zX/(1+|h|) \geq 1$ and $zX/(1+|h|) \leq 1$. For the last statement, we notice that the integral in question is 
\begin{displaymath}
  \frac{1}{|\Delta(4m-h^2)|^{1/4}}\int_0^{\infty} \tilde{g}(h; c) \frac{dc}{c}.
\end{displaymath}
By (2.4) and Lemma 8, we can remove the cut-off functions $\omega$ and replace $\tilde{g}$ by $g^{\ast}$ at the cost of 
\begin{displaymath}
 \ll \int_0^{(1+|h|)X/P} \frac{\min(1, \sqrt{c})}{X^{10}(1+|h|)^{5/2}} \frac{dc}{c} + \int_{(1+|h|)X^{10}}^{\infty} X \left(\frac{X}{c}\right)^{k-1} \frac{dc}{c} \ll \frac{1}{X^9(1+|h|)^{9/4}} 
\end{displaymath}
for $k \geq 4$. Now we interchange the order of integration in this absolutely convergent double integral. The inner integral over $c$ equals
\begin{displaymath}
  \int_0^{\infty} J_{k-1}\left(\frac{4\pi y\sqrt{m}}{c}\right) \cos\left(\frac{2\pi yh}{c}\right)\frac{dc}{c}
\end{displaymath}
which by \cite[6.693.2]{GR} is 0 for all $y > 0$ if $|h| > 2 \sqrt{m}$ and $k$ is even.\\ 

We are now prepared for the endgame. Obviously  $\phi_h$ satisfies the required properties for the application of Kuznetsov's formula for the congruence subgroup $\Gamma_0(4N')$. We use a smooth partition of unity to localize $z \asymp Z$ with $X^{-10} \ll Z \ll P/X^{1-\varepsilon}$. There are $\ll \log X$ of such terms. We call the truncated weight function $\phi_{h, Z}$. Now we apply Proposition 2; depending on the signs of $-\Delta$ and $4m-h^2$ as well as on the parity of $k$ we also use  (2.14) and/or (2.16). Precisely, let us first assume that $k$ and $\chi$ are even. If $-\Delta > 0$, we apply Proposition 2 with $\kappa = 3/2$; if $-\Delta < 0$, we use (2.14) first, and then apply Proposition 2 with $\kappa = 1/2$. If $k$ and $\chi$ are odd, we apply (2.16) first, and argue as above with $\kappa$ replaced by $2-\kappa$. In all cases, we obtain a sum over the spectrum of $\Delta_{\kappa}$ with $\kappa \in \{1/2, 3/2\}$. As we shall see, the exceptional eigenvalue $\lambda = 3/16$ (i.e.\ $t=i/4$) requires special care.  Let us first assume that $k$ and $\chi$ are even, and let us  distinguish four cases depending on the signs of $-\Delta$ and $4m-h^2$.  

a) If $-\Delta$ and $4m-h^2$ are both positive, then  by Proposition 2, the main term in (3.7) is a sum over
\begin{displaymath}
 \sum_{\substack{h \equiv s \, (2)\\ h \geq 0}} \left. \sum_j\right.^{(\frac{3}{2})}  \frac{4\sqrt{|\Delta(4m-h^2)|}\rho_j(|\Delta|)\bar{\rho}_j(4m-h^2) }{\cosh(\pi t_j)} \widehat{\phi}_{h, Z}(t_j) + \text{ two similar terms}
\end{displaymath}
(corresponding to the holomorphic and the continuous spectrum) for various values of $Z$. Let us write $\psi_Z(t) := \max_h |\widehat{\phi}_{h, Z}(t)|$. By Cauchy-Schwarz, this is at most
\begin{displaymath}
 \ll \sum_{|h| < 2\sqrt{m}} \left( \left.\sum_j\right.^{(\frac{3}{2})} \frac{|\Delta| |\rho_j(|\Delta|)|^2} {\cosh(\pi t_j)} \psi_Z(t_j)\right)^{1/2} \left( \left.\sum_j\right.^{(\frac{3}{2})} \frac{(4m-h^2) |\rho_j(4m-h^2)|^2}{ \cosh(\pi t_j)}   \psi_{Z}(t_j)\right)^{1/2}
\end{displaymath}
By the fourth statement of Lemma 6 with $Y := Z \leq 1$ and the second statement of Lemma 9, the exceptional eigenvalue $3/16$ (i.e.\ $t = i/4$) contributes $\ll Z^{1/2} \Xi(1) \ll (ZX)^{1/2} \ll P^{1/2}X^{\varepsilon}$; note that the $h$ and $j$ sum are bounded. Together with (2.8) and the second statement of Lemma 6 we see analogously that the other exceptional eigenvalues contribute $\ll Z^{-\theta}\Xi(1) \ll X^{1/2 + \theta+\varepsilon}$. Let us now turn towards the real $t_j$ It is customary to split the two $j$-sums into dyadic intervals $t_j \asymp T$. Now we use Lemma 5, the second statement of Lemma 6 (with $\theta = 0$ and $P := 1+XZ$, $Y := Z$) and Lemma 9. If $T \leq (1+XZ)X^{\varepsilon}$, we estimate each of the $j$-sums restricted to $t_j \asymp T$ by $\Xi(1) T^{3/2+\varepsilon}$ for bounded $h$. Choosing $A$ large enough in Lemma 6, we see that larger $T$  are negligible, so that we obtain a total contribution of
$\Xi(1) (1+ZX)^{3/2+\varepsilon} \ll P^{3/2}X^{1/2+\varepsilon}$. The same bound holds for the holomorphic and the continuous spectrum. 
 
b)  If $-\Delta > 0$ and $4m-h^2 < 0$,   the main term in (3.7) is of the shape 
\begin{displaymath}
 \sum_{\substack{h \equiv s \, (2)\\ |h| > 2\sqrt{m}}} \left. \sum_j \right.^{(\frac{3}{2})} \frac{4\sqrt{|\Delta(4m-h^2)|}\rho_j(|\Delta|)\bar{\rho}_j(4m-h^2) }{\cosh(\pi t_j)} \check{\phi}_{h, Z}(t_j) + \text{ continuous spectrum}.
\end{displaymath}
The exceptional eigenvalue $3/16$ does not contribute here, since all forms in this eigenspace have only positive Fourier coefficients. All other eigenvalues with $|t_j| \leq 2$, say,  contribute by (2.8), Lemmas 5, 6, and 9 (with $\min(A^{k-1}, A^{-1/2}) \leq A^{2+\varepsilon}$ for $A =zX/(1+|h|)$) at most 
\begin{displaymath}
   \sum_h |h|^{1/2} Z^{-\theta}\Xi(h) \ll \sum_h \frac{Z^{3/2-\theta+\varepsilon} X^{2+\varepsilon}}{(1+|h|)^{1+\varepsilon}} \ll P^{3/2-\theta+\varepsilon}X^{1/2+\theta + \varepsilon}.
 \end{displaymath}  
For the rest of the spectrum, we  cut the $h$-sum and the $j$-sum (resp.\ the $t$-integral) into dyadic pieces with $|h| \asymp H$, $t_j \asymp T$ and put $\psi_{Z, H}(t) := \max_{|h| \asymp H}|\check{\phi}_{h, Z}(t)|$. By Cauchy-Schwarz we get for each such subsum
\begin{equation}
\begin{split}
  \sum_{h \asymp H} & \left(\left. \sum_{|t_j| \asymp T}\right.^{(\frac{3}{2})} \frac{|\Delta| |\rho_j(|\Delta|)|^2 T^{3/2}}{\cosh(\pi t_j)}\psi_{Z, H}(t_j) \right)^{1/2}\\
& \times  \left( \left.\sum_{|t_j| \asymp T}\right.^{(\frac{3}{2})}\frac{ (h^2-4m) | \rho_j(4m-h^2)|^2 T^{-3/2}}{\cosh(\pi t_j)}  \psi_{Z, H}(t_j)\right)^{1/2}. 
\end{split}
\end{equation}
We estimate both factors using Lemmas 5, 7 and 11 getting
\begin{displaymath}
  \sum_{|h| \asymp H} \Xi(H) \left((1+ZX)^{3/2} \left((1+ZX)^{3/2} + H (1+ZX)^{1/2}\right)\right)^{1/2}. 
\end{displaymath}
We sum this over dyadic values for $H$. Splitting into $H \leq 1+XZ$ and $H > 1 +XZ$, we get a total contribution of $(1+ZX)^2 Z^{1/2}X^{1+\varepsilon} \ll P^{5/2}X^{1/2+\varepsilon}$. The same bound holds for the holomorphic and for the continuous spectrum. At this place there might be room for a small improvement by trying to treat the $h$-sum non-trivially.

  
c) The case $-\Delta < 0$ and $4m-h^2 > 0$ is analogous to the preceding case except that we apply (2.14) first and and the use  the Kuznetsov formula with weight $\kappa = 1/2$. In fact, this case is a little easier than the preceding one as the $h$-sum is bounded. 
 
d) For the last case, $-\Delta < 0$ and $4m-h^2 < 0$,  we use again (2.14) and  Proposition 2 with $\kappa = 1/2$. Thus we need to bound
\begin{displaymath}
 \sum_{\substack{h \equiv s \, (2)\\ h > 2 \sqrt{m}}}  \left. \sum_j\right.^{(\frac{1}{2})}  \frac{4\sqrt{|\Delta(4m-h^2)|}\rho_j(|\Delta|)\bar{\rho}_j(h^2-4m) }{\cosh(\pi t_j)} \widehat{\phi}_{h, Z}(t_j) + \text{ two similar terms}.
\end{displaymath}
We start with the contribution of the exceptional eigenvalue $\lambda_1= 3/16$ corresponding to $t_1=i/4$. To this end, we remove the partition of unity by summing over the various $Z$-pieces getting
\begin{displaymath}
\text{const }\cdot \sum_{\substack{h \equiv s \, (2)\\ h > 2 \sqrt{m}}}   \sqrt{h^2-4m}\bar{\rho}_1(h^2-4m)  \widehat{\phi}_{h}(i/4).
\end{displaymath}
By the last statement of  Lemma 6, Lemma  9 and our present assumption $h > 2 \sqrt{m}$ and $k$ even, we have
\begin{equation}
  \widehat{\phi}_h(i/4) \ll \frac{1}{X^9(1+|h|)^{5/2}} + \int_0^{PX^{\varepsilon-1}} \min(1, z^2) \Xi(h) \frac{dz}{z^{3/2}} \ll \frac{P^{5/2}}{(1+|h|)^{2} X^{1-\varepsilon}}.
\end{equation}
We cut the $h$-sum into pieces with $h \asymp H$; by Cauchy-Schwarz and Lemma 5, we have
\begin{displaymath}
  \sum_{h \asymp H}   \sqrt{h^2-4m}\bar{\rho}_1(h^2-4m) \ll H^{3/2}.
\end{displaymath}
Thus the total contribution of the eigenvalue $3/16$ is at most $P^{5/2}X^{\varepsilon-1}$. For the rest of the spectrum the analysis is exactly is in case b) above. For the other eigenvalues with $|t_j| \leq 2$ we use the bound $\widehat{\phi}_{h, Z}(t_j) \ll Z^{-\theta} \Xi(h)$ and obtain similarly a total contribution of $P^{3/2-\theta+\varepsilon}X^{1/2+\theta+ \varepsilon}$. For the remaining parts of the spectrum we cut once again the $h$-sum and the $j$-sum  into pieces $|h| \asymp H$, $|t_j| \asymp T$,  and put $\psi_Z(t) = \max_{|h| \asymp H} \widehat{\phi}_{h, Z}(t)$. As in (3.8), we  need to bound
\begin{displaymath}
\begin{split}
  \sum_{h \asymp H} & \left(\left. \sum_{|t_j| \asymp T}\right.^{(\frac{1}{2})} \frac{|\Delta| |\rho_j(|\Delta|)|^2 }{\cosh(\pi t_j)}\psi_{Z, H}(t_j) \right)^{1/2}
   \left( \left.\sum_{|t_j| \asymp T}\right.^{(\frac{1}{2})}\frac{ (h^2-4m) | \rho_j(h^2-4m)|^2}{\cosh(\pi t_j)}  \psi_{Z, H}(t_j)\right)^{1/2}. 
\end{split}
\end{displaymath}
which by Lemmas 5, 6 and 9 is at most
$P^{5/2}X^{1/2+\varepsilon}$. The same bound holds for the holomorphic and for the continuous spectrum. This completes the proof if $k$ is even; in particular $c=0$ in this case.\\ 

If $k$ and $\chi$ are odd, we apply (2.16) first and argue similarly. The cases b) and c) do not cause any difficulty, in case d) the eigenvalue $\lambda=3/16$ is now negligible by the same argument is in case a) above. In particular, $c=0$ if $\Delta > 0$. Finally case a) is as above except that now there might be a large contribution of the exceptional eigenvalue $\lambda = 3/16$. Precisely, let $u$ run through an $L^2$-normalized basis $\mathcal{B}$ of $\textbf{H}_{\text{exc}} := \textbf{H}_{1/2}(4N', \chi\chi_{-4}, i/4)$. We remove the partition of unity, that is, we re-sum over dyadic values of $Z$.  Our present assumption is $-\Delta > 0$ and $4m-h^2> 0$, so by (2.12) and (3.6), the eigenvalue $\lambda=3/16$ contributes
\begin{displaymath}
 \sum_{\substack{h \equiv s \, (2)\\ 0 \leq h < 2 \sqrt{m}}}  \sum_{u \in \mathcal{B}}  4\sqrt{2|\Delta|(4m-h^2)}\rho_u(|\Delta|)\bar{\rho}_u(4m-h^2)\widehat{\phi}_{h}(i/4).
\end{displaymath}
By Lemmas 6 and 9 we see similarly as in (3.9) that this is 
\begin{displaymath}
  (1+i)\sum_{\substack{h \equiv s \, (2)\\ 0 \leq h < 2 \sqrt{m}}}  \sum_{u \in \mathcal{B}}  4\sqrt{2|\Delta|(4m-h^2)}\rho_u(|\Delta|)\bar{\rho}_u(4m-h^2) \int_0^{\infty} \phi_h(x) x^{-3/2} dx + O(P^{5/2}X^{\varepsilon-1}).
\end{displaymath}
We insert the definition (3.6) of $\phi_h$ getting
\begin{displaymath}
 4i\sum_{\substack{h \equiv s \, (2)\\ 0 \leq h < 2 \sqrt{m}}}  \sum_{u \in \mathcal{B}}  4(4|\Delta|(4m-h^2))^{1/4}\rho_u(|\Delta|)\bar{\rho}_u(4m-h^2) \int_0^{\infty} \tilde{g}\left(h; \frac{\sqrt{|\Delta|(4m-h^2)}}{4x}\right)  \frac{dx}{x}.
\end{displaymath}
By Lemma 8 we can remove the cut-off functions $\omega$ and replace $\tilde{g}$ by $g^{\ast}$ at  the cost of a negligible error. Substituting (3.3) we arrive at a contribution of
\begin{equation}
\begin{split}
&4i\sum_{\substack{h \equiv s \, (2)\\ 0 \leq h < 2 \sqrt{m}}}  \delta(h) \sum_{u \in \mathcal{B}}  4(4|\Delta|(4m-h^2))^{1/4}\rho_u(|\Delta|)\bar{\rho}_u(4m-h^2)\\
&\times \int_0^{\infty} \int_0^{\infty} w\left(y-\frac{s}{2}\right) J_{k-1}\left(\frac{16 \pi xy\sqrt{m}}{\sqrt{|\Delta|(4m-h^2)}}\right) \cos\left(\frac{8\pi hxy}{\sqrt{|\Delta|(4m-h^2)}} \right)    \frac{dxdy}{x}.
\end{split}
\end{equation}
By \cite[6.693.2]{GR} the double integral equals
\begin{equation}
\int_0^{\infty}w(y)dy \int_0^{\infty} J_{k-1}(2 x \sqrt{m})\cos(xh) \frac{dx}{x} = \frac{1}{k-1}T_{k-1}\left(\frac{h}{2\sqrt{m}}\right) \int_0^{\infty}w(y)dy 
\end{equation}
where $T_{\nu}(x) = \cos(\nu \arcsin(x))$ is the Chebychev polynomial. Recall that $\delta(h)$ was defined as 2 if $h \not=0$ and $1$ otherwise. The preceding two displays give an expression for the constant $c$ in Theorem 2 if $k$ is odd, $\Delta < 0$ and\footnote{Recall that we have found an explicit description of $c$ in (3.1) if $\Delta = 0$.}  if $f = P_m, m \not= \square$.  As explained at the beginning of the section, every cusp form can be written as a linear combination of these $P_m$ up to a remainder that contributes at most $O(1)$ to (1.3).  The constant  agrees with (1.4) in the introduction: Observe that in this case the Whittaker function satisfies $W_{1/4, -1/4}(4 \pi ny) = (4\pi n y)^{1/4} e^{-2 \pi n y}$ \cite[(4.21)]{DFI}, so the Fourier coefficients $\rho_u(n)$ and $a_j(n)$ are related by $a_j(n) = (4\pi n)^{1/4} \rho_u(n)$. This completes the proof of Theorem 2.\\

Using results of \cite{SS}, the Fourier coefficients of $u \in \textbf{H}_{\text{exc}} = \{y^{1/4}f \mid f \in M_{1/2}(4N', \chi\chi_{-4})\}$ can be described explicitly. In \cite{SS}, an explicit basis of $f \in M_{1/2}(4N', \chi\chi_{-4})$ is constructed in terms of theta-functions of the type $\sum_n \psi(n) e(bn^2 z)$ where $\psi$ the (even) primitive Dirichlet character underlying $\chi\chi_{-4b}$ of conductor $r$, say, and $r^2b \mid N'$. This gives further conditions on the vanishing of the constant $c$. For example,  we see that $\rho_u(|\Delta|) = 0$ for all $u \in \textbf{H}_{\text{exc}}$ unless $|\Delta| = -\Delta = b \square$ for some $b \mid N'$. Moreover, if $N'$ is not divisible by an odd square (other than 1) and not by 256, then $r \mid 8$, hence $\psi$ is real, and $\textbf{H}_{\text{exc}} \not= \{0\}$ unless $\chi$ is real.  Let us conclude with the specific (and easy to generalize) 

\textbf{Example.} Let $f = P_3 \in S_5(12, \chi_{-4})$ as in (2.1) and $q(x) = x^2 +x+1$. Then $\Delta = -3$, $N=N' = 12$, $m=3$ and $k = 5$, so $T_4(x) = 8x^4-8x^2+1$.  By the above discussion, $M_{1/2}(48, \textbf{1})$ is generated by certain $\sum_n \psi(n) e(bn^2 z)$ where $\psi$ has conductor $r$ and $r^2b \mid 12$. This implies $r \mid 2$, hence $r=1$, so $\psi$ is trivial. Moreover, we need $bn^2 = 3$ for some $n$, otherwise $\rho_u(|\Delta|) = 0$. Thus only the normalized version of $\theta(z) := \sum_{n \in \Bbb{Z}} e(3n^2z)$ contributes to (3.10), and hence only $h=3$ contributes to the sum.  Let 
\begin{displaymath}
  \alpha = \| \theta \|^2 = \int_{\Gamma_0(48)\backslash \mathcal{H}} |\theta(z)|^2 y^{1/2} \frac{dxdy}{y^2}.
\end{displaymath}
Then the constant $c = c(f, q)$ in Theorem 2 is given by
\begin{displaymath}
  c = 4i \cdot 2 \cdot 4(4\cdot 3 \cdot 3 )^{1/4} \frac{1}{\alpha} \frac{1}{4} T_4\left(\frac{\sqrt{3}}{2}\right) = -\frac{4\sqrt{6}i}{\alpha}.
\end{displaymath}
\\

Theorem 1 is now immediate. In fact, it is enough to show 
\begin{displaymath}
  \sum_{X/2 \leq n \leq X} a(q(n))  = cX/2 +O_{f, q, \varepsilon}(X^{6/7+\varepsilon})
\end{displaymath}
We approximate the characteristic function on $[X/2, X]$ by a function $w$ that is 1 on $[X/2 + X/P, X-X/P]$ and 0 on $[0, X/2] \cup [X, \infty)$ and satisfies $w^{(j)} \ll_j (P/X)^j$ for all $j \in \Bbb{N}_0$. By Deligne's bound (that holds for non-newforms as well) this introduces an error of $O(X^{1+\varepsilon}/P)$. We estimate the smoothed sum by Theorem 2 and equalize the error terms by choosing $P= X^{1/7}$. This completes the proof of Theorem 1. 

\section{Proof of Proposition 2}

Let us finally give the postponed proof of the second formula in Proposition 2. For $m \in \Bbb{N}$ and $\Re s > 1$ define non-holomorphic Poincar\'e series by
\begin{displaymath}
  P_m(z, s):=  \sum_{\gamma \in \Gamma_{\infty} \backslash \Gamma_0(N)}\bar{\chi}(\gamma)  j(\gamma, z)^{-2\kappa} (\Im \gamma z)^s e(m\gamma z). 
\end{displaymath}
Its Fourier expansion is given by (see e.g.\ \cite[(15)]{Pr})
\begin{displaymath}
  P_m(z, s) = y^se(mz)+y^s\sum_{\ell \in \Bbb{Z}} e(\ell x) \sum_{N \mid c} \frac{K_{\chi, \kappa}(m,\ell; c)}{c^{2s}} B(c, m, \ell, y, s)
\end{displaymath}
where
\begin{displaymath}
  B(c, m, \ell, y, s) = \int_{\Bbb{R}} e\left(-\frac{m}{c^2z} - \ell x\right) \left(\frac{z}{|z|}\right)^{-\kappa} \frac{dx}{|z|^{2s}}
\end{displaymath}
where $z = x + iy$. The formula in question will be proved by calculating the inner product $\langle U_m(., s_1), \overline{U_{n}(., s_2)}\rangle$ in two ways: Using the Fourier expansion and unfolding the fundamental domain, one finds as in  \cite[Lemma 4.3]{DI} that\footnote{Note that there is a factor 2 in the first display of p.\ 252 of \cite{DI} that is missing in the statement of Lemma 4.3.}
\begin{equation}
\begin{split}
  \langle U_m(., s_1), \overline{U_{n}(., s_2)}\rangle = &\frac{\pi 2^{3-s_1-s_2}e(-\frac{\kappa}{4}) \Gamma(s_1+s_2 -1)}{\Gamma(s_1-\frac{\kappa}{2})\Gamma(s_2+\frac{\kappa}{2})}\left(\frac{m}{n}\right)^{\frac{s_2-s_1}{2}}\\
& \times\sum_{N \mid c} \frac{K_{\chi, \kappa}(m, -n; c)}{c^{s_1+s_2}}K_{s_1-s_2}\left(\frac{4\pi\sqrt{mn}}{c}\right)
\end{split}
\end{equation}
for $\Re s_1, \Re s_2 > 1$. For convenience we sketch the argument: By the unfolding technique, we get
\begin{displaymath}
  \langle U_m(., s_1), \overline{U_{n}(., s_2)}\rangle = \sum_{N \mid c} \frac{K_{\chi, \kappa}(m, -n; c)}{c^{2s_1}} \int_0^{\infty} B(c, m, -n, y, s_1)y^{s_1+s_2-2}e^{-2\pi ny}dy. 
\end{displaymath}
We substitute the integral for $B$, reverse the order of integration and change variables $x = y\xi$ getting
\begin{displaymath}
  \int_{-\infty}^{\infty} (1+\xi^2)^{-s_1}\left(\frac{i+\xi}{|i+\xi|}\right)^{-\kappa} \int_0^{\infty} y^{s_2-s_1-1} e\left(-\frac{m}{c^2y(i+\xi)}+(i+\xi)ny\right)dyd\xi.
\end{displaymath}
The inner integral equals (cf.\ \cite[3.471.10]{GR})
\begin{displaymath}
  2 \left(\frac{m}{nc^2}\right)^{(s_2-s_1)/2} (1-\xi)^{s_1-s_2} K_{s_1-s_2}\left(\frac{4\pi \sqrt{mn}}{c}\right).
\end{displaymath}
Now the $\xi$-integral can be expressed in terms of $\Gamma$-functions yielding (4.1). By Weil's bound (2.15) this holds still in $\Re s_1, \Re s_2 \geq 5/6$, say. On the other hand, using the spectral theorem, one shows as in \cite[Lemma 2]{Pr}
\begin{displaymath}
\begin{split}
 & \langle U_m(., s_1), \overline{U_{n}(., s_2)}\rangle =   \frac{(4\pi)^{2-s_1-s_2}m^{1-s_1}n^{1-s_2}}{\Gamma(s_1-\frac{\kappa}{2})\Gamma(s_2+\frac{\kappa}{2})}\left(\left.\sum_j\right.^{(\kappa)} \bar{\rho}_j(m)\rho_j(-n) \prod_{\pm, \, i=1, 2} \Gamma(s_i-\frac{1}{2}\pm it_j) \right. \\
  & \left.+ \frac{1}{4\sqrt{mn}} \sum_{\mathfrak{a} \text{ singular}}\int_{\Bbb{R}} \left(\frac{n}{m}\right)^{it} \frac{\bar{\phi}_{\mathfrak{a}}(m, \frac{1}{2}+it) \phi_{\mathfrak{a}}(-n, \frac{1}{2}+it) }{\Gamma(\frac{1+\kappa}{2} - it)\Gamma(\frac{1-\kappa}{2} + it)} \prod_{\pm, \, i=1, 2} \Gamma(s_i-\frac{1}{2}\pm it)dt\right).
\end{split}
\end{displaymath}
Equating the last two expressions for $\langle U_m(., s_1), \overline{U_{n}(., s_2)}\rangle$ with $s_1 = 1+ir$ and $s_2 = 1-ir$, $r \in \Bbb{R}$, we get
\begin{equation}
\begin{split}
   e\left(-\frac{\kappa}{4}\right) & \sum_{N \mid c}\frac{K_{\chi, \kappa}(m, -n; c)}{c} x K_{2ir}(x)\\
& = 4\sqrt{mn}\left.\sum_j\right.^{(\kappa)} \frac{\bar{\rho}_j(m)\rho_j(-n)}{\cosh(\pi t_j)} \frac{\pi^2 \cosh(\pi t_j)}{2\cosh(\pi(r+t_j))\cosh(\pi(r-t_j))}  \\
   & + \sum_{\mathfrak{a} \text{ singular}}\int_{\Bbb{R}} \left(\frac{n}{m}\right)^{it} \frac{\bar{\phi}_{\mathfrak{a}}(m, \frac{1}{2}+it) \phi_{\mathfrak{a}}(-n, \frac{1}{2}+it) }{\cosh(\pi t ) \Gamma(\frac{1+\kappa}{2} - it)\Gamma(\frac{1-\kappa}{2} + it)}  \frac{\pi^2\cosh(\pi t)}{2\cosh(\pi(r+t))\cosh(\pi(r-t))} dt
\end{split}
\end{equation}
where we have set $x := 4\pi \sqrt{mn}/c$. To conclude the proof of the proposition, we appeal to the Kontorovich-Lebedev formula \cite[p.\ 361]{Sn}: If $f$ is a smooth function on $\Bbb{R}_{>0}$ such that $xf(x)$ and $x (x^{-1}f(x))'$ are absolutely integrable, 
and if we put 
\begin{displaymath}
  L_f(r) = \int_0^{\infty} K_{ir}(y) f(y)\frac{dy}{y},
\end{displaymath}
then $f$ can be recovered by 
\begin{displaymath}
  f(x) = \frac{1}{\pi^2}\int_{-\infty}^{\infty} L_f(r) K_{ir}(x) \sinh(\pi r)  r dr.
\end{displaymath}
For $\phi$ as in Proposition 2, we multiply (4.2) by
\begin{displaymath}
  \frac{4}{\pi^2} r \sinh(2\pi r) \int_0^{\infty} K_{2ir}(y)\frac{\phi(y)}{y} \frac{dy}{y}
\end{displaymath}
and integrate over $r \in (-\infty, \infty)$. By the Kontorovitch-Lebedev inversion formula with $f(x) = \phi(x)/x$, the left hand side becomes
\begin{displaymath}
  e\left(-\frac{k}{4}\right) \sum_{N \mid c}\frac{K_{\chi, \kappa}(m, -n; c)}{c} \phi(x).
\end{displaymath}
For the right hand side we use the inversion formula for $f(x) = K_{2ir}(x)x$ together with the formula (a special case of \cite[6.576.4]{GR})
\begin{displaymath}
  \int_0^{\infty} K_{2ir}(y)K_{2it}(y) dy = \frac{\pi^2}{4\cosh(\pi(r-t))\cosh(\pi(r+t))}
\end{displaymath}
which completes the proof of the second formula of Proposition 2.\\

\textbf{Remark:} This formula holds for arbitrary weight $\kappa \in [0, 2)$ if one uses Kloosterman sums with an appropriate multiplier. It is interesting to note that the opposite sign formula is much less sensitive to different weights than the same sign formula.


\begin{thebibliography}{BHM2}

\bibitem[Bi]{Bi} A. Bir\'o, \emph{Cycle integrals of Maass forms of weight 0 and Fourier coefficients of Maass forms of weight 1/2}, Acta Arith. \textbf{94} (2000), 103-152

\bibitem[Bl]{Bl} V. Blomer, \emph{On the central value of symmetric square $L$-functions}, to appear in Math. Z. 


\bibitem[By]{By} V. A. Bykovski\u i, \emph{Spectral decompositions of certain automorphic functions and their number-theoretic applications}, J. Sov. Math. \textbf{36} (1987), 8-21

\bibitem[DI]{DI} J.-M. Deshouillers and H. Iwaniec, \emph{Kloosterman sums and Fourier coefficients of cusp forms}, Invent. math. \textbf{70} (1982), 219-288


\bibitem[DFI]{DFI} W. Duke, J. Friedlander and H. Iwaniec, \emph{The subconvexity problem for Artin $L$-functions}, Invent. math. \textbf{149} (2002), 489-507

\bibitem[FI]{FI} J. Friedlander and H. Iwaniec, \emph{A polynomial divisor problem}, J. Reine Angew. Math. \textbf{601} (2006), 109-137


\bibitem[GR]{GR} I. S. Gradshteyn, I. M. Ryzhik, \emph{Tables of integrals, series, and products}, 5th edition, Academic Press, New York, 1994


\bibitem[Ho]{Ho} C. Hooley, \emph{On the number of divisors of quadratic polynomials}, Acta. Math. \textbf{110} (1963), 97-114

\bibitem[Iw]{Iw} H. Iwaniec, \emph{Fourier coefficients of modular form of half-integral weight}, Invent. math. \textbf{87} (1987), 385-401
 
\bibitem[IK]{IK} H. Iwaniec, E. Kowalski, \emph{Analytic number theory}, American Mathematical Society Colloquium Publications 53, American Mathematical Society, Providence, RI, 2004

\bibitem[KS]{KS} H. Kim, {\em Functoriality for the exterior square of $GL\sb 4$ and the symmetric fourth of $GL\sb 2$ (with Appendix 1 by D. Ramakrishnan and Appendix 2 by H. Kim and P. Sarnak)}, J. Amer. Math. Soc. \textbf{16} (2003), 139--183

\bibitem[Ki]{Ki} H. Kim, \emph{Functoriality and number of solutions of congruences}, Acta Arith. \textbf{128} (2007), 235-243


\bibitem[Mo]{Mo} Y. Motohashi, \emph{Spectral theory of the Riemann zeta function}, Cambridge 1997

\bibitem[Pr]{Pr} N. V. Proskurin, \emph{On the general Kloosterman sums}, J. Math. Sci. (New York) \textbf{129} (2005), 3874-3889


\bibitem[SS]{SS} J.-P. Serre and H. M. Stark, \emph{Modular forms of weight 1/2}, Modular functions of one variable VI, Lecture notes in Mathematics 627 (1977), 27-67 

\bibitem[Sh1]{Sh} G. Shimura, \emph{On modular forms of half-integral weight}, Ann. of Math \textbf{97} (1973), 440-481

\bibitem[Sh2]{Sh2} G. Shimura, \emph{On the holomorphy of certain Dirichlet series}, Proc. London Math. Soc. (3) \textbf{31} (1975), 79-98

\bibitem[Sn]{Sn} I. N. Sneddon, \emph{The use of integral transforms}, McGraw Hill, New York 1972

\bibitem[T]{T} V. Tipu, \emph{Polynomial Divisor Problems}, Ph.D. thesis, University of Toronto 2008
  
\end{thebibliography}
\end{document}